# A Solution to a problem and the Diophantine equation $x^2 + bx + c = y^2$

*Konstantine Hermes Zelator*
*Department of Mathematics*
*College of Arts and Sciences*
*Mail Stop 942*
*University Of Toledo*
*Toledo, OH 43606-3390*
*U.S.A.*

# 1 Introduction

In the September 2007 issue of the *Newsletter of the European Mathematical Society* *(*see **[1]** *)*, the following problem appeared, whose solution was solicited : Does the diophantine equation $x^2 + x + 1 = y^2$ have any nontrivial integer solution for $x$ and $y$ ? A solution is offered in **[2]**. That solution, an easy one indeed, can be obtained as follows:

Write $x^2 + x + 1 - y^2 = 0$, as a quadratic in $x$. Then, $x = \dfrac{-1 \pm m}{2}$, for some non-negative integer $m$ such that $m^2 = 4y^2 - 3$. Note that, in order for $x$ to be an integer, it is necessary that $m$ is rational.

However, $m = \sqrt{4y^2 - 3}$, the square root of an integer. The square root of an integer is either a positive integer; or otherwise a positive irrational number. More generally, the $n^{th}$ root of a positive integer is either a positive integer or an irrational number. Equivalently, an integer is equal to the $n^{th}$ power of a rational number if and only if, it is the $n^{th}$ power of an integer. For more details, the interested reader can refer to [3] or [4]. Now, with $m$ being an integer, we obtain

$$(m - 2y)(m + 2y) = -3$$

Obviously, the possibilities are rather limited. Either one of the two factors above is equal to -1 and other to 3; or alternatively one of them must equal 1, the other -3. This couched with the condition $m \geq 0$, and a straightforward calculation yields $y = 1$ and $m = 1$; or $y = -1$ and $m = 1$. Accordingly, we obtain the values $x = 0, -1$. The four pairs $(0,1), (0,-1), (-1,1), (-1,-1)$; constitute all the solutions of the above Diophantine equation. Two of them are trivial ( in the sense $xy = 0$), the other two nontrivial.

Note that $-3$ is the discriminant of the quadratic trinomial $x^2 + x + 1$. This problem provided the motivation for tackling the more general problem: Finding all the integer solutions; that is, all pairs $(x, y)$ in $\mathbf{Z} \times \mathbf{Z}$, which satisfy the diophantine equation,

$$x^2 + bx + c = y^2 \qquad (1),$$

for given integers b and c.

The approach we employ, is the same; as in the case b = c=1 above. We only use basic knowledge and material found in the early part of a first course in elementary number theory, therefore accessible to second or third year university students with a major (or a minor) in mathematics.

## 2 The diophantine Equation $x^2 + bx + c = y^2$

Equation (1) is equivalent to,

$$x^2 + bx + c - y^2 = 0 \qquad (2)$$

If we consider (2) as a quadratic equation in $x$; it has discriminant

$D(y) = b^2 - 4c + 4y^2$

As in the case $b = c = 1$, $D(y)$ must be an integer square; $D(y) = m^2$, for some $m \in \mathbf{Z}$, with $m \geq 0$. Accordingly, (2) gives,

$$x = \frac{-b \pm m}{2}, m \geq 0, b^2 - 4c = (m - 2y)(m + 2y) \qquad (3)$$

By (3), $b^2 \equiv m^2 \pmod 4$, which implies $b \equiv m \pmod 2$.

Note that we must have $b \equiv m \pmod 2$, since $x$ is an integer.

## 3 Analysis

*A. Assume $b^2 - 4c \neq 0$ in (3).*

According to (3), to find all the integral solutions of (2), we must proceed as follows: For each pair $(\delta_1, \delta_2)$ of divisors $\delta_1$ and $\delta_2$ of the integer $b^2 - 4c$; and with $\delta_1 \delta_2 = b^2 - 4c$,

$$\delta_1 \delta_2 = b^2 - 4c,$$

We simply take $m + 2y = \delta_1$ and $m - 2y = \delta_2$. So that,

$$\left. \begin{array}{l} \delta_1 \delta_2 = b^2 - 4c, 0 \leq m = \dfrac{\delta_1 + \delta_2}{2}, y = \dfrac{\delta_1 - \delta_2}{4}) \\ \\ \text{and with } \delta_1 \equiv \delta_2 \pmod 4 \end{array} \right\} \qquad (4)$$

A straightforward calculation shows that if an integer pair $(x, y)$, satisfies (3) and (4); then it must be a solution of (1). If $b$ is odd, then $b^2 \equiv 1 \pmod 4$ and so $b^2 - 4ac \equiv 1 \pmod 4$. In this case(i.e., $b$ odd), it is easily seen in (4) that we must have either $\delta_1 \equiv \delta_2 \equiv 1 \pmod 4$; or alternatively, $\delta_1 \equiv \delta_2 \equiv 3 \pmod 4$.

On the other hand, if $b$ is even, the situation is somewhat more complicated. Indeed, put $b = 2B$, for some integer $B$. Then, $b^2 - 4c = 4(B^2 - c)$. If $B^2 - c$ is an odd integer; then all integer solutions to (2), can be obtained by taking $\delta_1 = 2d_1, \delta_2 = 2d_2$ where $d_1$ and $d_2$ are odd integers such that $d_1 d_2 = B^2 - c$. In this manner, all integral solutions of

(2) can be produced. Note that from (4) and (3), we obtain

$m = d_1 + d_2 \geq 0, x = B \pm (d_1 + d_2)$, and $y = \dfrac{d_1 - d_2}{2}$.

If $B^2 - c$ is an even integer, then this case splits into two subcases. If $B^2 - c \equiv 2 \pmod 4$, there exist no integer solutions to (2); since then $\delta_1 \delta_2 = 4(B^2 - c) \equiv 8 \pmod{16}$; which implies that either one of $\delta_1, \delta_2$ is a multiple of 4; and the other is congruent to $2 \pmod 4$. Or, alternatively one of $(\delta_1, \delta_2)$ is a multiple of 8( more precisely, exactly divisible by 8), while the other one is odd. In either case, the condition $\delta_1 \equiv \delta_2 \pmod 4$, fails to be satisfied as required by (4). Next, if $B^2 - c \equiv 0 \pmod 4$, the solutions of (2) are readily obtained by taking $\delta_1 = 4d_1, \delta_2 = 4d_2$, and with $d_1 d_2 = \dfrac{B^2 - c}{4}$. In this subcase, we have $0 \leq m = 2(d_1 + d_2), x = -B \pm 2(d_1 + d_2)$ and $y = d_1 - d_2$. Observe that the number of integral solutions of (2) is always even. This is clear when $b^2 - 4c$ is not equal to minus an integer or perfect square; because then $m > 0$ ( since $m \geq 0; m = 0$ or $m > 0$. But $m = 0$ implies by (4), that $m$ is minus an integer square), and so for each integer value of $y$ in (4), there are two distinct integer values for $x$; hence the number of solutions must be even. Now, when $b^2 - 4c = -(\text{integer square}) = -k^2$; for some integer k>0. Then, the number of integral solutions of (2) is equal to $N_1 + N_2$ ; where $N_1$ is the number of solutions with $m > 0$; and $N_2$ the number of solutions with $m = 0$. Obviously $N_1$ is even ; for it is either zero, or if it is positive; it must be even for the same reason given above for the case when $b^2 - 4c$ is not equal to minus a perfect square. To find $N_2$, note that if (2) has solutions when $m=0$; we have $\delta_1 = -\delta_2$ in (4). And so $b^2 - 4c = -\delta^2 = -k^2$ which yields,

$\delta_1 = k$ and $\delta_2 = -k$; or $\delta_1 = -k$ and $\delta_2 = k$. Two solutions $\left(-\dfrac{b}{2}, \dfrac{k}{2}\right)$ and $\left(-\dfrac{b}{2}, -\dfrac{k}{2}\right)$.

(Note that $b$ must be even ; and so must $k$). So we see that $N_2 = 2$.

*B. Assume $b^2 - 4c = 0$ in (3)*

If $b^2 - 4c = 0; c = \dfrac{b^2}{4}$, which requires $b$ to be even. And so (1) $\Leftrightarrow \left(x + \dfrac{b}{2}\right)^2 = y^2$; and thus $y = \pm\left(x + \dfrac{b}{2}\right)$. Thus, we see that when $b^2 - 4c = 0$; the diophantine equation has infinitely many solutions.

# 4  Conclusions

As a result of the previous section, we can state the following theorem.

**Theorem 1**

Consider the diophantine equation $x^2 + bx + c = y^2$

1. Assume $b^2 - 4c \neq 0$. Then this equation has an even number of solution pairs in
   **Z** x **Z**  and,

(i)    If $b$ is an odd integer, then all the solutions of the above equation are given by

$$x = \dfrac{-2b \pm (\delta_1 + \delta_2)}{4}, y = \dfrac{\delta_1 - \delta_2}{4};$$

where $\delta_1, \delta_2$ can be any integers such that $\delta_1 \delta_2 = b^2 - 4c, \delta_1 + \delta_2 > 0$ and with either $\delta_1 \equiv \delta_2 \equiv 1 \pmod 4$; or $\delta_1 \equiv \delta_2 \equiv 3 \pmod 4$.

Note that $-2b \equiv 2 \equiv \pm(\delta_1 + \delta_2)(\bmod 4) \Rightarrow -2b \pm (\delta_1 + \delta_2) \equiv 0(\bmod 4)$.

(ii) If $b$ is even; $b = 2B$ and $B^2 - c$ is an odd integer; then all integer solutions are given by,

$$x = -B \pm (d_1 + d_2), y = \frac{d_1 + d_2}{2}$$

where $d_1, d_2$ can be any odd integers such that $d_1 d_2 = B^2 - c$ and $d_1 + d_2 > 0$.

(iii) If $b$ is even; $b = 2B$, and $B^2 - c \equiv 2(\bmod 4)$; the above equation has no integer solutions.

(iv) If $b$ is even; $b = 2B$ and $B^2 - c \equiv 0(\bmod 4)$ all integer solutions are given by,

$$x = -B \pm 2(d_1 + d_2), y = d_1 - d_2$$

where $d_1$ and $d_2$ can be any integers such that $d_1 d_2 = \frac{B^2 - c}{4}$ and $d_1 + d_2 \geq 0$.

2. If $b^2 - 4c = 0$, the above equation has infinitely many solutions given by

$$y = \pm\left(x + \frac{b}{2}\right).$$

# 5 A Proposition

## Proposition 1

Let $p$ be an odd prime and suppose that the integers $b$ and $c$ satisfy $b^2 - 4c = p^{n-1}$, for some positive integer $n$. Then the diophantine equation $x^2 + bx + c = y^2$ has exactly $2n$ solutions.

*Remark 1:* It is easy to generate the numbers $b, c, p$ of proposition 1. Indeed, pick $b$ odd. If $n$ is even, then $n-1$ is odd, and let $p$ be a prime, $p \equiv 1 \pmod 4$. Then $b^2 \equiv 1 \equiv p^{n-1} \pmod 4$. If $n$ is odd, then $n-1$ is even; then take $p$ odd prime; $b^2 \equiv 1 \equiv p^{n-1} \pmod 4$. In either case, take $\dfrac{b^2 - p^{n-1}}{4} = c$.

**Proof :** Since $p$ is odd; $p^{2n}$ is odd, and thus, $b$ must be odd as well. We apply Theorem 1, part(i). From $\delta_1 \delta_2 = b^2 - 4c = p^{n-1}$ and $\delta_1 + \delta_2 > 0$; we see that $\delta_1$ and $\delta_2$ must both be positive. Since $b^2 - 4c$ is a prime power we easily see that the choices for the divisors $\delta_1$ and $\delta_2$ are $\delta_1 = 1, p, p^2, \ldots, p^{n-1}$.

And correspondingly $\delta_2 = p^{n-1}, p^{n-2}, p^{n-3}, \ldots, 1$.

In summary, $(\delta_1, \delta_2) = (p^i, p^{n-1-i})$; for $i = 0, 1, \ldots, n-1$.

Obviously, for each of the $n$ pairs $(\delta_1, \delta_2)$, exactly one integer value of $y$ is obtained; while two distinct integer values for $x$. "It is straightforward to show that the $n$ $y$-values obtained in the manner, are distinct (left to the reader)". A total of $2n$ solution pairs $(x, y)$. It follows from proposition 1, that for each positive integer $2n$, there exist infinitely many diophantine equations of the form of (1); which have exactly $2n$ solutions. Consequently we have the following

**Corollary 1 :** *For each non-negative even integer $2n$, there exist infinitely many diophantine equations of the form $x^2 + bx + c = y^2$ with exactly $2n$ solutions.*

# 6 Exactly two solutions

A cursory look at the formulas in (3) and (4) reveals that, equation (1) will have exactly two solution pairs in one of the two ways:

Either both pairs have the same $x$-value and opposite $y$-values; or both pairs have distinct $x$-values and the same $y$-value. In the first case, $m = 0$ and so $b^2 - 4c = -k^2$, for some positive integer $k$. In the second case, $b^2 - 4c = k^2$. In the either case, $k$ cannot contain an odd prime factor; for in such a case, equation (1) would have atleast four solutions. Thus, $k$ must have a power of 2. A bit more work establishes that (1) will have exactly two solutions precisely when $b^2 - 4c = 1, 4, -4, 16,$ or $-16$. We provide the following tabulation:

| | | |
|---|---|---|
| 1 | $b^2 - 4c = 1$ | Solutions are $(x, y) = \left(\dfrac{1-b}{2}, 0\right), \left(-\dfrac{(b+1)}{2}, 0\right)$ |
| 2 | $b^2 - 4c = 4$ | Solutions are $(x, y) = \left(\dfrac{2-b}{2}, 0\right), \left(-\dfrac{(b+2)}{2}, 0\right)$ |
| 3 | $b^2 - 4c = 16$ | Solutions are $(x, y) = \left(\dfrac{4-b}{2}, 0\right), \left(-\dfrac{(b+4)}{2}, 0\right)$ |
| 4 | $b^2 - 4c = -4$ | Solutions are $(x, y) = \left(-\dfrac{b}{2}, -2\right), \left(-\dfrac{b}{2}, 2\right)$ |
| 5 | $b^2 - 4c = -16$ | Solutions are $(x, y) = \left(-\dfrac{b}{2}, -4\right), \left(-\dfrac{b}{2}, 4\right)$ |

# 7 Closing Remarks

(a) Equation (1) has always infinitely many rational solutions; for any integer values of $b, c$. Indeed, if $r_1, r_2$ are any rationals such that $r_1 r_2 = b^2 - 4c \neq 0$ and $m = \frac{r_1 + r_2}{2} \geq 0$.

Then, $(x, y) = \left( \frac{-b \pm m}{2}, \frac{r_1 - r_2}{4} \right)$ is a rational solution of (1); all the rational solutions can be found in this manner. If $b^2 - 4c = 0$; $(x, y) = \left( x, \pm \left( x + \frac{b}{2} \right) \right); x \in Q$ are all the rational solutions.

(b) If we allow $b, c$ to vary over the integers, we see that equations (1) describe a family of hyperbolas; each hyperbola having center $\left( -\frac{b}{2}, 0 \right)$, and axes of symmetry the lines $x = -\frac{b}{2}$ and $y = 0$; and with asymptotes the lines $y = \pm \left( x + \frac{b}{2} \right)$. Each such hyperbola has an even number of integral points (per Theorem 1), but an infinite number of rational points.

(c) Some hyperbolas contain in fact, no rational points. For example, the curves $x^2 + 1 = 3y^2$, $y^2 + 1 = 3x^2$. More generally, one can prove that if $p$ is a prime, $p \equiv 3 \pmod 4$; the curve $x^2 + 1 = py^2$ has no rational points (we invite the reader to show this).

Hint: First show that if a rational point on the curve $x^2 + 1 = py^2$ exists; then this implies $m^2 + n^2 = pk^2$; for some positive integers $m, n, k$. Then use the fact that -1 is a quadratic non residue of $p$; to prove that this such $m, n, k$ cannot exist.

(d) Some hyperbolas, on the other hand, contain an infinite number of integral points. For example, $x^2 - 2y^2 = 1$ is one of them (well known Pell-equation). More generally, $x^2 - dy^2 = 1$; where $d$ is a positive integer, not a perfect square.

# 8 References

[1] EMS Newsletter, September 2007, *Problem 18, page 46*.

[2] Ovidiu Furdui and Konstantine Zelator, *A Solution to Problem 18 (EMS Newsletter, September 2007)*.

[3] Kenneth H. Rosen, *Elementary Number Theory and Its Applications*, third edition, 1993, Addison-Wesley Publishing Co., ISBN: 0-201-57889-1. *See page 96, Theorem 2.11*.

[4] W. Sierpinski, *Elementary Theory of Numbers*, Warsaw, 1964. ISBN: 0-598-52758-3. *See page16, Theorem 7, and the Corollary*.